\documentclass{amsproc}

\newcommand{\RR}{{\mathbb R}}
\newcommand{\NN}{{\mathbb N}}
\newcommand{\ZZ}{{\mathbb Z}}
\newcommand{\TT}{{\mathbb T}}
\newcommand{\CC}{{\mathbb C}}
\newcommand{\EE }{{\mathbb E}}

\newcommand{\cI}{{\mathcal I}}
\newcommand{\cL}{{\mathcal L}}

\newcommand{\cO}{{\mathcal O}}

\newcommand{\cT}{{\mathcal T}}

\newcommand{\ua}{\underline a}

\newcommand{\up}{\underline p}

\newcommand{\ur}{\underline r}
\newcommand{\us}{\underline s}
\newcommand{\ut}{\underline t}
\newcommand{\uu}{\underline u}
\newcommand{\uv}{\underline v}
\newcommand{\ux}{\underline x}
\newcommand{\uxi}{\underline \xi}

\newcommand{\uA}{\underline A}

\newcommand{\ru}{\mathrm{u}}
\newcommand{\rv}{\mathrm{v}}
\newcommand{\rh}{\mathrm{h}}
\newcommand{\rf}{\mathrm{f}}
\newcommand{\rA}{\mathrm{A}}
\newcommand{\rB}{\mathrm{B}}
\newcommand{\rF}{\mathrm{F}}

\newcommand{\Id}{\mathrm{Id}}

\newcommand{\D}{\partial }
\newcommand{\eps}{\varepsilon } 
\newcommand{\ra}{\rangle}
\newcommand{\la}{\langle} 
\newcommand{\VV}{\vert\!\vert\!\vert} 

\newtheorem{theorem}{Theorem}[section]
\newtheorem{lemma}[theorem]{Lemma}
\newtheorem{prop}[theorem]{Proposition}
\newtheorem{cor}[theorem]{Corollary}

\theoremstyle{definition}
\newtheorem{definition}[theorem]{Definition}
\newtheorem{example}[theorem]{Example}

\theoremstyle{remark}
\newtheorem{remark}[theorem]{Remark}

\newtheorem{assumption}[theorem]{Assumption}

\numberwithin{equation}{section}

\begin{document}

\title{ Remarks on  the Well-Posedness of 
the Nonlinear Cauchy Problem }

\author{Guy M\'etivier}
\address{MAB, Universit\'e de Bordeaux I, 33405 Talence Cedex, France}
 \email{metivier@math.u-bordeaux.fr}
\thanks{Research partially   supported by European network HYKE,  HPRN-CT-2002-00282 .}

 
\subjclass{Primary 35 }
\date{January 1, 1994 and, in revised form, June 22, 1994.}

\dedicatory{This paper is dedicated to Fran\c{c}ois Treves.}

\keywords{Cauchy Problem, well posedness}

\begin{abstract}
We show that hyperbolicity is  a necessary condition for
  the well posedness of the noncharacteristic Cauchy problem
 for nonlinear partial differential equations. 
 We 
 give  conditions on the initial data which are necessary for the  existence of solutions
 and we analyze Hadamard's instabilities in Sobolev spaces.  
 We also show that genuinely nonlinear
 equations raise new interesting  problems. 
\end{abstract}

\maketitle

\section{Introduction}

The question of the well-posedness of the Cauchy problem was first raised by  Hardamard 
 (\cite{Ha}, \cite{Ha1}) who proved that it is ill-posed in the case of linear second order elliptic  equations. But  the introduction  in \cite{Ha1}  clearly indicates that 
Hadamard was interested in   nonlinear equations as well. 
In modern words, Hadamard's proof  is based on the analytic regularity 
of linear elliptic boundary problems. This   
regularity  has been extended to nonlinear elliptic equations  by  Morrey  (\cite{MN}) 
so that    Hadamard's argument 
also applies  to general nonlinear elliptic equations.  

For general {\sl linear} equations, 
it is well known that   hyperbolicity is a necessary condition 
for the well-posedness of the noncharacteristic  Cauchy problem  in $C^\infty$,   
that is for the existence of solutions for general $C^\infty$ data  
(see Lax \cite{La}, Mizohata \cite{Mi} and    Ivrii-Petkov \cite{IP}  for 
a simplified proof and further developments;   see also \cite{Ho}). 
Moreover, for several classes of nonhyperbolic equations, explicit necessary conditions on the initial data  for the existence of solutions have been given (see  
  \cite{Ni1} \cite{Ni2}). 
For {\sl  nonlinear} equations, Wakabayashi \cite{Wak}  has proved 
that the existence of a smooth {\it stable} solution implies hyperbolicity,  
  stability meaning that one can perturb the initial data and 
the  source terms in  the 
 equations. In a previous paper, Yagdjian obtained  this result, with  a much weaker definition of stability, 
 in the sense of continuous dependence on the initial data, for the special case of ``gauge invariant''
  equations (\cite{Yag1}). 
We also mention  \cite{Yag} for a particular case and \cite{Hou} for 
first order scalar complex equations.

In this paper, we continue the analysis of  Hadamard's instabilities
 for nonhyperbolic nonlinear equations in two directions. 
 First, in the spirit of Hadamard's examples, we give necessary conditions on the 
 initial data for the existence of smooth solutions of a given equation, without perturbing
 the equation or the source terms. 
 Next, we also want to point out that 
the nonlinear  theory yields interesting 
and difficult new problems. 
There are many interesting examples, for instance in  multi-phase fluid dynamics, 
where  the  equations  are not {\sl everywhere} hyperbolic. 
To mention one occurrence of this phenomenon, consider  Euler's equations of gas dynamics in Lagrangian coordinates:
\begin{equation}
\label{ex11}
\left\{\begin{aligned}
&\partial_t u  +  \partial_x v \, =\, 0\, , 
\cr
& \partial_t v + \, \partial_x p(u) \, =\, 0\, .   
\end{aligned}\right.
\end{equation}
The system  is hyperbolic  [resp. elliptic]  when
$p'(u) > 0$ [resp. $p'(u) < 0$].  
For van der Waals state laws, it happens that $p$   is decreasing on an interval
$[u_*, u^*]$.  
A mathematical example is 
\begin{equation}
p(u) = u \, (u^2 - 1)\, 
\label{ex12}
\end{equation}
Hadamard's argument (see e.g. \cite{Ha}, \cite{Ha1}) shows that the Cauchy problem with data taking values in the elliptic region is ill-posed: if  $u_{\vert t = 0} $ is real analytic near $\ux$
and $u(\ux)$ belongs to  the elliptic interval, then 
 any local 
$C^1$  solution is analytic (see e.g. \cite{MN}); thus the initial data $v_{ \vert t = 0} $ {\sl must} be analytic for 
the initial value problem to have a solution. 
  
 However,  there are other interesting questions about the system \eqref{ex11}. 
There are  classical solutions 
with values in one of the hyperbolic region
$u < u_*$ or $u > u^*$, 
but there are  also discontinuous solutions, 
modeling for instance phase transitions,   which  take values in both regions. 
They  have been extensively studied, see e.g. \cite{Sl}, \cite{Be}, \cite{Fr}.  
  Another  remarkable fact is that \eqref{ex11} has a conserved 
    energy.  Let $P(u)$ satisfy $P' = p$. 
    Then  the energy: 
   \begin{equation}
E(t) \, := \,   \int   \big({1\over 2}\,  v^2(t,x)  + P(u (t,x))\big) \, dx \,   
\end{equation}
 is  conserved for solutions of \eqref{ex11}. 
    For the example   \eqref{ex12}, 
$  P(u)  = {1 \over 4} \, u^2 \, (u^2 - 2)$. If one considers the 
    periodic problem,   the $L^2$ norm is dominated by
   the $L^4$ norm on $[0, 2 \pi]$, thus the  boundedness of $E$ 
    controls the $L^4$ and $L^2$ norm of the solutions.
    Of course, this is formal, and  the validity of  {\sl a priori} bounds does not   prove
    the {\sl existence} of solutions. However, this indicates that the nonexistence of solutions
    is  much more subtle than in the linear case. In particular,  there is no blow up phenomenon 
    in $L^p$ norms.

   The equations \eqref{ex11} are thought of as approximations or limits of
   more complicated models which may include for instance viscosity or capillarity
   (see e.g. \cite{Be}); 
    numerical schemes have also been considered (see e.g. \cite{Fan}, \cite{HW}).  
    In the case of periodic solutions, 
  spectral methods lead to filter high frequencies and to consider 
the ``approximate'' system
\begin{equation}
\label{ex14}
\left\{\begin{aligned}
&\partial_t u^\lambda   +  \partial_x v^\lambda   =  0  , 
\\
& \partial_t v^\lambda  +  \partial_x S_\lambda  p(u^\lambda )   = 0,    
\end{aligned}\right. \qquad
\left\{\begin{aligned}
&  u^\lambda{}_{\vert t = 0} = S_\lambda h , 
\\
& v^\lambda{}_{\vert t = 0} = S_\lambda k,  
\end{aligned}\right.
\end{equation}
where $S_\lambda$ is the projector on Fourier modes of index 
$\vert n \vert \le \lambda$.   For instance, when $p$ is given by \eqref{ex12}, 
the conservation of energy and the Cauchy-Lipschitz theorem  imply that
 in the periodic case: 
   \begin{enumerate}
  \item[] {\sl  for all $h \in L^4$ and $k \in L^2$, the equations 
  $\eqref{ex14}$ have  global solutions  $(u^\lambda, v^\lambda)$ which 
  are uniformly bounded in $C^0([0, \infty[ ; L^4 \times L^2)$.} 
  \end{enumerate}
Note that there are no conditions on $h$, which can take values in the elliptic region 
   $u^2 < 1/3$. The question is to analyze the behavior of 
   $(u^\lambda, v^\lambda)$ as $\lambda \to + \infty$. 
   Because of the bounds on  $(u^\lambda, v^\lambda)$ and 
    $(\D_t u^\lambda, \D_t v^\lambda)$, there are subsequences 
    which converge weakly and strongly in $C^0([0, 1], H^{- \eps})$. In particular,
      the weak limits $(\uu, \uv)$ are bounded with values in $L^2$ and continuous
    in time for the weak tolopogy of $L^2$. Thus  $\uu(0) = h$ and 
    $\uv (0) = k$. Taking the weak limit  
    $\up$ of $p(u^\lambda)$,  there holds
   $$
   \D_t \uu + \D_x \uv = 0 \,, \quad \D_t \uv + \D_x \up = 0\,. 
   $$
The  question is to express  $\up$ in terms of $\uu$ and $\uv$.  
  As mentioned above,  the  answer,  
   $\up = p(\uu)$,   $(\uu, \uv)$ smooth,  
   {\sl cannot} be true in general  when  $h$ takes values in the elliptic zone.  
The common idea is that the  limits $(\uu, \uv)$  ``escape'' from the elliptic region,
as suggested by numerical  calculations  (\cite{Fan}, \cite{HW}), but 
no  rigorous proof of this fact seems available in the literature.  
A detailed answer to the questions above remains  a very interesting 
open problem. Motivated by this problem, we consider  in Section 5 
 a modified nonlocal system: 
 \begin{equation}
\label{ex15}
\left\{\begin{aligned}
&\partial_t u  =  a(t)  \partial_x v  , 
\cr
& \partial_t v = \vert a(t) \vert  \partial_x u  ,  
\end{aligned}\right. \qquad  \mathrm{with}\quad 
a(t) = \Vert u(t) \Vert^2_{L^2}  -1 . 
\end{equation}
 This a  version of Kirchhoff equations (\cite{K}, \cite{Li}
 \cite{AS}), which is non hyperbolic when 
  $a  < 0$.  As \eqref{ex11}, this system has a natural (formal) energy:   
 \begin{equation}
 E (t) = \Vert v(t, \cdot)  \Vert^2_{L^2} + \Big\vert  \Vert u(t, \cdot)  \Vert^2_{L^2} - 1  \Big \vert . 
 \end{equation}
 This implies that the equations with filtered initial data $(S_\lambda h, S_\lambda k)$ 
 has global  solutions  $(u^\lambda, v^\lambda)$ uniformly bounded in 
 $L^\infty([0, + \infty[, L^2(\TT))$. 
For   large classes of  ``nonanalytic'' initial data $(h, k) \in L^2 $,    
with $\Vert h \Vert_{L^2} < 1$, we show that 
the limits are $\uu = h$, $\uv = k$, constant in time, remaining  in the elliptic region. 
The limit equations, $\D_t \uu = \D_t \uv = 0$, have little to see with the original ones. 
This indicates that the answers to the questions above might be very delicate. 
 
\bigbreak 
 
  Now we review the results of  Sections 
 2 to 4. 
 To fix a framework, we consider  first order square systems
\begin{equation}
\partial_t u  = F(t, x, u, \partial_x u) \, , \quad u_{\vert t = 0} = h \, . 
\label{CP}
 \end{equation} 
 where 
 $F $ is a smooth function of 
$(t, x, u, v) \in \RR \times \RR^d \times \RR^N \times (\RR^N)^d$.  
The principal  symbol 
of the equation reads
\begin{equation}  
\tau \Id  - \xi \cdot \partial_v F (t, x, u, v) \, . 
\end{equation}
Hyperbolicity means  that all the eigenvalues of 
$\partial_v F$   are real. 
We consider the local Cauchy problem \eqref{CP} near 
$(0, \ux)$ and a given  base point  $(\uu, \uv)$, assuming that the initial data 
satisfy 
\begin{equation} 
h(\ux)  =  \uu \, , \quad \partial_x h(\ux )  =  \uv . 
\label{D0}
\end{equation}
 The results in  Sections 2 to 4  illustrate   
 the idea that
if $\D_v F(0, \ux, \uu, \uv)$ has a  non real eigenvalue, then 
the Cauchy problem \eqref{CP} \eqref{D0} 
for classical solutions  is ill-posed.

 Well posedness means first {\sl solvability}.   
  Hadamard's  counterexamples (see   the example \eqref{ex11} above)  
  prove that analyticity type conditions on the data are necessary for the existence of solutions
  of the elliptic Cauchy problem.   
 In the same vein, 
 consider the  equation\footnote{
 This explicit elementary example was suggested by  
 Nicolas Lerner.}:
\begin{equation}
  (\D_x+ i \D_y ) u = u^2 \,, \quad  x > 0,  \qquad  u_{\vert x= 0} = h\,, 
\end{equation} 
with $h(0) \ne 0$.  Any $C^0$ solution on $\{ x \ge 0 \} $  near the origin is $ C^1$, 
does not vanish and   $1/u$ satisfies the   equation 
$ 
(\D_x+ i \D_y )(1/u) = -1/2\,. 
$ 
Therefore, $1/u + \overline z/2 $ is holomorphic and $1/h $ is the trace of an holomorphic function in $\{ x > 0 \}$, implying that 
necessarily  $1/h$ is microlocally analytic in the direction $+ 1$   at the origin. 
In particular, if $h$ is real valued, $h$ must be real analytic near the origin.

 In   Section 2,  we extend  this analysis 
 to   first order scalar complex  equations. 
 For such equations, is is proved in \cite{Hou}  that the existence of solutions 
 for all complex data close to a given $h$, implies that the system
 must be semilinear and hyperbolic.  
 We give a more precise result, showing that,  in the non hyperbolic case,   
{\sl microlocal analyticity conditions on the initial data are necessary for the existence of classical solutions},  
hence that the Cauchy problem has no classical solution for most  initial data. 
The proof is based on the analysis of \cite{Me1} (see also \cite{BGT}), which provides 
approximate integral representation of the $C^1$ or $C^2$ solutions. 
Taking real and imaginary parts, this provides examples of $2 \times 2$ systems
where the Cauchy problem has no classical solutions. 

This analysis  does not extend to   general systems:  the representation
and approximation  theorems valid in  the scalar case have no analogue;  
the local uniqueness of the Cauchy problem may be false (\cite{Me2});  
  there are no microlocal analytic regularity  theorem at elliptic directions for 
$C^1$ or $C^2$ solutions. 
However, there are results about $H^{s'}$ microlocal regularity for 
$H^{s}$ solutions when $s' \le 2s - s_0$ (see Bony \cite{Bo} and   Sabl\'e-Tougeron \cite{ST}).
In Section 4,  we show   that 
 if $\D_v F(0, \ux, \uu, \uv )$ has a nonreal eigenvalue, then for 
all $H^s $ local solution, with $s > \frac{d}{2} + 1$, 
 the polarized $H^{s'}$ wave front set of the initial data 
 is not arbitrary, when $s' < 2s - \frac{d}{2} - 1$. 
 In particular,   for most data in $H^{s'}$, the Cauchy problem has 
 no local  $H^s$ solution. 
 Note that $s'$ can be taken larger  than $s$, 
and   for any  ``loss''  $k $, it applies to   $s' = s + k$, if $s$ is large enough.  
The restriction $s > \frac{d}{2} + 1$ is natural in order to have a $C^1$ classical solution. 
The restriction $s' < 2s - d -1/2$  is forced  by the Theorem of 
microlocal ellipticity for nonlinear equations (see \cite{Bo} \cite{ST}).
Under additional assumptions on the equations, one can 
show that for arbitrary large  $s'$ there are $H^{s'} $ initial data such that 
the Cauchy problem has no local  $H^{s}$ solution.

\medbreak

For linear equations,  standard functional analysis results convert  
 well-posedness into estimates, and necessary conditions are found 
 by contradicting the estimates.  Solvability implies 
continuous dependence on the data (see  also F.John \cite{Jo} for general 
remarks on this notion). For nonlinear equations, 
there are no such abstract argument and it is reasonable to include the continuous 
dependence in the definition of the well posedness.  
In addition, because  local   uniqueness is not 
guaranteed, we also include  it in  the following  definition of H\"older continuous 
solvability. 
In the next statement, 
$B_r$ denotes the ball $\{ \vert x - \ux \vert < r \}$ 
and $\Omega_{r, \delta} $ the lens shaped domain 
$$
\Omega_{r, \delta}  = \big\{ (t, x) \ : \  0 <    t \,,    \vert x -\ux  \vert^2 + \delta t < r^2 
\big\}
$$

\begin{definition}
\label{def12}
We say that the Cauchy problem $\eqref{CP}$ is H\"older well posed 
on  $H^\sigma$, if there are constants $r_0 > r_1 > 0$,  $\delta> 0$, $c > 0$, $C$ and   
$\alpha \in ]0, 1]$, such that for all 
$h \in H^\sigma (B_{r_0})$  satisfying 
$\Vert h - \uu - x \cdot \uv \Vert _{H^\sigma(B_{r_0})} \le c$ and 
all $r \in ]0, r_1]$, the Cauchy problem $\eqref{CP}$ has 
a unique solution in $C^1(\overline \Omega_{r, \delta})$, 
with norm bounded by $C$. Morever, given $h_1$ and $h_2$, the 
corresponding solutions satisfy for all $r < r_1$: 
\begin{equation} 
\Vert u_1 - u_2 \Vert_{L^\infty(\Omega_{r, \delta}) } 
\le  C \Vert h_1 - h_2 \Vert^\alpha_{H^\sigma(B_{r_0})} 
\label{Hest}
\end{equation}
\end{definition}

In Section 3, we show that {\sl if $F$ is real analytic
and  if $\D_v F(0, x_0, u_0, v_0)$ has a nonreal eigenvalue, then 
the Cauchy problem $\eqref{CP}$ is not H\"older well posed 
on $H^\sigma$, for all $\sigma$}. 

Note that the definition above differs strongly from 
the notion of {\it stable} solution introduced in \cite{Wak} in the sense that 
we do not allow perturbations of the equations, while 
the stability used in \cite{Wak}  is related to the solvability of 
\begin{equation}
\D_t u + F(t, x, u,   \D_x u ) = f \,, \quad  u_{\vert t = 0} = h + h'.  
\label{CPc}
\end{equation}
for all $f$ and $h'$ small. 
The analysis  is based  on the construction of 
 asymptotic solutions  using WKB or geometric  optics expansions. 
 But they are not exact solutions, yielding error terms $f$ 
which  are precisely the  source terms 
chosen in \cite{Wak}.  In this analysis, the choice of 
$f$ is dictated by the choice of $h$. 
It is  interesting and much stronger to consider exact solutions of \eqref{CP}
(as in \cite{Yag1}), or to be able to choose 
$h$ and $f$ independently. 
 In Section 3, we construct exact solutions close to the  approximate solutions,
 by Cauchy-Kowalewski type arguments. This is where we use the analyticity of the equation. 
 In this respect, the results of this section  give a detailed account of  the  $H^s$ instability
 of analytic solutions when hyperbolicity fails. 


\section{Necessary conditions   for scalar complex equations} 

 To simplify the discussion, consider a quasilinear scalar equation 
\begin{equation}
\partial_t u \, +  \, \sum_{j=1}^d  a_j(t, x, u) \partial_{x_j} u + b(t, x, u) \, = \, 0\, , \quad
u_{\vert t=  0} = h\, .
\label{eq21} 
\end{equation}
where the  $a_j $ are holomorphic functions of $(t,x,u)$ on a neighborhood
of $(0, \ux, \uu )$. The Cauchy data $h$ is always assumed to satisfy
$h(\ux) = \uu$.     
The nonhyperbolicity condition reads
\begin{equation}
  {\rm Im}\,  a(0, \ux , \uu  ) \ne 0\, . 
\label{nH}
\end{equation} 
\bigbreak 

\begin{theorem} 
\label{th21}
If the Cauchy problem $\eqref{eq21}$ has a  
$C^1$ solution for $t \ge 0$ on a neighborhood of  $(0, \ux)$, then  for all 
$\xi \in \RR^d  $ such that $\xi \cdot   {\rm Im}\,  a(0, \ux , \uu  )  > 0$, 
$(\ux, \xi) $
does not belong to the analytic wave front set of  $h$. 
\end{theorem}
\bigbreak 
 
For the definition of the analytic wave front set, we refer to \cite{Sj} or 
\cite{Ho}.  In particular,  it contains the $C^\infty $ wave front set
and the theorems implies that if the local Cauchy problem
has a $C^1$ solution,  then $h$ must be $C^\infty $ at 
$(\ux, \uxi)$ if $\uxi \cdot   {\rm Im}\,  a(0, \ux , \uu  )  > 0$. 
This means that for all 
$C^\infty$ cut-off function $\chi$ supported in
a sufficiently   small neighborhood of $\ux$, 
the Fourier transform of $\chi h$ is rapidly decreasing 
in any small conical neighborhood of $\uxi$. 
For ``most''   functions $h$ in $H^s$, $(\ux, \uxi)$  belongs to the 
 the $C^\infty $ wave front set. 
Theorem \ref{th21} implies that for  most $h$,  the Cauchy problem 
\eqref{eq21} has no $C^1$ solution.

\begin{example} 
\textup{Taking real and imaginary parts of the unknowns yields nonexistence theorem
for $2 \times 2$ real systems. With 
$\alpha_j (u, v) = {\rm Re\, } a_j (u+ i v)$, $\beta_j (u, v) = {\rm Im\, } a_j(u+iv)$, 
the equation \eqref{eq21}  with $b = 0$ is equivalent to:}
\begin{equation}
\left\{\begin{aligned}
& \D_t u +  \sum \D_{x_j} \alpha_j(u, v)  = 0 \,,  \quad  u_{\vert t = 0} = h,
\\
& \D_t v +   \sum \D_{x_j} \beta_j(u, v) = 0 \,, \quad v_{\vert t = 0} = k. 
\end{aligned}\right. 
\label{ex22}
\end{equation}
\textup{Suppose that $\underline \beta = \beta(h(\ux), k(\ux) ) \ne 0$ and 
choose $\uxi$ such that $\uxi\cdot \underline \beta > 0$.  If 
$h+ ik$ is not microlocally analytic at $(\ux, \underline \xi )$, then the Cauchy problem  
\eqref{ex22} has no local $C^1$ solution near $(0, \ux)$. }

\textup{ For instance, this applies to the the system:}
\begin{equation}
\left\{\begin{aligned}
& \D_t u +   u \D_x u - v \D_x v +  \D_y u  = 0 \,,  \quad  u_{\vert t = 0} = h,
\\
& \D_t v +   v \D_x u + u \D_x v + \D_y v  = 0 \,, \quad v_{\vert t = 0} = k, 
\end{aligned}\right. 
\label{ex23}
\end{equation}
\textup{when $k(\ux ) \ne 0$. 
For functions independent of $y$, or equivalently dropping
the $\D_y$,  the system is elliptic for $v \ne 0$ and Hadamard's argument applies. 
The example \eqref{ex23} shows that Theorem \ref{th21} also applies  to  nonelliptic systems. }
\end{example}
\bigbreak

\begin{proof}[Proof of Theorem $\ref{th21}$]

{\bf a)}    
The complex characteristic curves are  integral curves of  
the holomorphic vector field: 
$$
 \mathcal{ L} = \partial_t  + \sum a_j(t,x,u) \partial_{x_j} \, -  \, 
b(t,x,u) \partial_u\, . 
$$
They are given by  
$Z_j (t, x, u) = c_j , U(t, x, u)  = c_0$, where  $Z_j$ and  $U$ are local 
holomorphic solutions of  
$$
\left\{ \begin{aligned} \mathcal{L}  Z_j = 0\, ,& \quad Z_j{}_{\vert t = 0} = x_j\, ,
\\
\mathcal{L}  U = 0\, ,& \quad U{}_{\vert t = 0} = u\, .
\end{aligned}\right.
$$ 
We also introduce the additional
variables  $v = (v_1, \ldots, v_d)$, which are placeholders for $\D_{x_j} u$, 
and the function 
$$
J(t,x,u,v) := \det \Big({ \partial Z_j (t,x,u) \over \partial{x_k} } + 
v_k {\partial Z_j(t,x,u) \over \partial u} \Big).
$$
Let  $G(t, x, u)$ be a  holomorphic solution of $\cL G = 0$ 
on a complex neighborhood $\cO$ of $(0, \ux, \uu)$. 
Suppose that $u$ is   $C^1$  solution  
 of \eqref{eq21} on  
$[0, T] \times \Omega$   such that $(t, x, u(t, x)) \in \cO$ 
for all $(t, x) \in [0,T] \times \Omega$. Then, 
by  Lemma 2.2.2 of  \cite{Me1},  there holds for all $s \in [0, T]$
and  
$\chi \in C^1_0(\Omega)$: 
\begin{equation}
\begin{aligned}
\int_\Omega & G(0, x, h(x) )  \chi(x) \, dx  = 
\int_\Omega G(s, x, u(s,x)) \, \chi(x)\,  \widetilde J (t,x) \, dx  \,   
\\ \, 
&  -  \int_{[0,s] \times \Omega}
\sum_j \partial_{x_j} \chi (x)\, a_j(t,x,u(t,x))\, G(t, x, u(t,x) ) 
\, \widetilde J (t, x) \,  dt dx
\end{aligned}
\label{eq23}
\end{equation}
with  $\widetilde J(t, x) := J(t, x, u(t, x),  \partial_xu (t, x))$.

\medbreak 
{\bf b)}
We use \eqref{eq23} with  
\begin{equation}
G_{\lambda, y}(t,x,u) := \Big({\lambda \over \pi}\Big)^{d/2}\, 
U(t,x,u)  \, e^{ - \lambda q ( Z(t,x,u) - y)}\, .
\label{eq24}
\end{equation}
where $q(y) = \langle Q y, y \rangle$ is a quadratic form, with real coefficients, 
positive definite on $\RR^d$. 
 
The $G_{\lambda, y}$  are defined and holomorphic for 
$\vert t \vert \le T$, $ \vert x - \ux  \vert \le r$, 
$\vert u- \uu \vert \le \rho $, for some 
$T> 0$, $ r > 0$  and $\rho  > 0$. 
We can also assume that the  given solution  $u$   of \eqref{eq21} 
is defined and $C^1$ for real  $(t, x) \in [0, T] \times \overline \Omega$
where $\Omega$ is the ball $ \{ \vert x - \ux \vert  <  r \}$ and that 
$\vert u(t, x) - \uu \vert < \rho $ on this domain. 
We fix   $\chi \in C^\infty_0(\Omega)$ equal to 1 on a smaller  
neighborhood of $\ux$. 
Because 
$
Z(t, x, u) = x + O\big(\vert t \vert  )\big) 
$, 
${\rm Re} Z (t, x, u(t, x) ) \ne 0$ for 
$t$ small and $x $ in the support of $d\chi$. 
Because $Z(0, x , u) - \ux  = x - \ux \ne 0$ on  the support of 
$d\chi$, 
there are  
$\Omega_0 \subset \Omega $, real neighborhood of     
$  \ux$,   $\varepsilon > 0$, 
$\delta > 0$ and $T_0 > 0$, such that 
\begin{equation}
\begin{aligned}
\forall y \in \Omega_0 +i [-\delta, \delta]^d, \ & \forall t  \in [0, T_0] , 
\  \forall x \in \mathrm{supp} d \chi \ : 
\\ \quad 
& {\rm Re}q  (Z(t, x, u(t, x)) - y) \ge 2 \varepsilon > 0. 
\end{aligned}
\label{eq25}
\end{equation}
  
  Consider 
\begin{equation}
T  h(y, \lambda) \, := \,  \Big({\lambda \over \pi}\Big)^{d/2}\, 
\int_\Omega e^{ -\lambda q(x-y) }\, h(x) \, \chi(x) \, dx.  
\label{eq26}
  \end{equation}
 We apply \eqref{eq23}  to $G= G_{\lambda, y}$ given by \eqref{eq24}.
  The  estimate \eqref{eq25} shows that the second integral 
  in the right hand side is 
  $O(e^{- \varepsilon \lambda})$.  Therefore, there is $C$ such that 
for all 
 $y \in \Omega_0 +i [-\delta, \delta]$ and  $t \in [0, T_0]$,
\begin{equation}
\begin{aligned}
\Big\vert T h (y, \lambda) -  
\Big({\lambda \over \pi}\Big)^{d/2}\, \int_\Omega 
\widetilde U(t,x)   e^{ - \lambda q (\widetilde Z(t,x) - y) } \chi(x)  
 \widetilde J (t, x)    &  dx \Big\vert 
\\  &  \le C e^{ -\varepsilon \lambda} \,  
\label{eq27}
\end{aligned}
\end{equation}
where   $\widetilde Z(t, x) := Z(t,x, u(t,x)) $  and we use similar notations for $\widetilde U$ and $\widetilde J$  
(see the estimate (4.3.1) in \cite{Me1}).

\medbreak 
{\bf c)} We now make use of Assumption \eqref{nH}. Shrinking 
$\Omega$ if necessary, in addition to the previous requirements, we can further assume that
\begin{equation}
\forall x \in \Omega \ : \quad 
\vert {\rm Im} (a (0, x, h(x) ) - \ua  \vert \le  \rho ,  
\end{equation}
where $\ua = {\rm Im} (a (0, \ux, h(\ux) )$ and 
  $\rho > 0$ to be chosen later on.  We use the estimate \eqref{eq27} with 
\begin{equation}
y = \ux -  i t {\rm Im} \ua  + y'\,,    \quad y' \in \CC^d,  \quad
 \vert y' \vert \le  \rho t .   
\label{eq29}
\end{equation}  
 Because   
$  Z(t,x, u)  = x -  t a (0,x,u) + O(t^2) $ and $u \in C^1([0, T] \times \overline \Omega)$, 
there holds:  
$$
{\rm Im } \big( \widetilde Z(t,x) - y \big) = - { \rm Im} y' -   t \big({\rm Im} a(0, x, h(x)) - \ua \big) +    O(t^2)\,. 
$$
Thus, there is $T_1 > 0$ such that for $t \in [0, T_1]$: 
\begin{equation*}
\vert {\rm Im }\big(  \widetilde Z(t,x) - y \big) \vert \le   3 \rho t .
\label{eq210}
\end{equation*}
Hence
\begin{equation*}
 q \big( {\rm Im }\big(  \widetilde Z(t,x) - y \big) \big) \vert \le   9 \Vert Q \Vert  \rho^2 t^2 .
\label{eq210}
\end{equation*}
On the other hand: 
$  {\rm Im} y   =  t \ua  -  {\rm Im }y' $ and therefore if $\rho$ is small enough 
and $\vert y' \vert \le \rho t$, 
\begin{equation*}
q ( {\rm Im} y ) \ge t^2 q (\ua) / 2 . 
\end{equation*}
 Hence, if 
 $\rho$ small enough, for  $t \in ]0, T_1]$,  $y$ satisfying  \eqref{eq29} and $x \in \Omega$ there holds:
\begin{equation}
 - {\rm Re} q (\widetilde Z(t, x) - y ) 
 \le    q \big( {\rm Im} (\widetilde Z(t, x) - y) \big)  \le  q \big( {\rm Im} y \big) 
 - \frac{t^2}{4}q( \ua )
 \label{eq211}
 \end{equation}
 We now fix $t > 0$, $t \le \min(T_0, T_1)$,  such that 
 $y \in \Omega_0 + i [- \delta, \delta]^d$ for all 
 $y$ satisfying \eqref{eq29}.   
 Thus, the estimates \eqref{eq27} and \eqref{eq211}   imply  that  there are 
 $\varepsilon_1 > 0$  and $C >0$  such that for all 
 $y$ in the complex ball of radius $\rho t   $ centered at 
 $\ux - it \ua$ and for all $\lambda \ge 1$,  there holds 
 \begin{equation}
 \label{n213}
 \vert T h(y, \lambda) \vert  \le C e^{ \lambda(  q ( {\rm Im} y )   - \varepsilon_1)   } \, .
\end{equation}
 Since the quadratic form $q$ is definite positive on $\RR^d$, 
 for $y \in \CC^d$, the unique real critical point of
 $x \mapsto \mathrm{Re} q(y - x) $  is $x = \mathrm{Re} y$ and at this point 
 $ - \D_x    q( y - x)$ is equal to $- Q \mathrm{Im} y$. 
By Proposition 7.2 of Sj\"ostrand \cite{Sj} (see also \cite{Del}, section I.2), 
 the estimate \eqref{n213}  on a neighborhood of $\ux - it \ua$   implies that 
 $(\ux, t Q \ua)$ does not belong to the analytic wave front set of $h$.  
 
 \medbreak
{\bf d) }  For all $\xi$ such that $  \xi \cdot \ua > 0$, there is a definite positive real symmetric  $Q$ such that  $Q \ua = \xi$. We apply the previous step to  
$q(x) = \langle Q x, x \rangle$  which implies that 
 there is $t > 0$ such that  $(\ux, t \xi)$ does 
not belong to the analytic wave front set of $h$. Since the wave front is conic in 
$\xi$, the theorem is proved. 
\end{proof}


\section{Hadamard's instabilities in Sobolev spaces}

We consider systems, and for simplicity we state the results for quasi-linear 
systems:
\begin{equation}
\D_t u = \sum_{j=1}^d A_j(t, x, u) \D_{x_j}  u + F(t,x, u) , \quad 
u_{\vert t = 0} = h. 
\label{eq31}
\end{equation}
We assume that the $A_j$ and $F$ are real valued and real analytic near 
$(0, \ux, \uu) \in \RR \times \RR^d \times \RR^N$. 
We want to compare two solutions of \eqref{eq31} with initial data $h_1$ and 
$h_2$. We can choose $h_1$ to be analytic, for instance 
$h_1 (x) = \uu$, and find an analytic  local solution 
$u_1$ by Cauchy-Kowalewski theorem. 
Changing $u$ to $u- u_1$, we get an equation similar 
to \eqref{eq31}, with the additional information that   $0$ is a solution, that is: 
\begin{equation} 
F(t, x, 0) = 0   \quad \mathrm{ or }  \quad F(t, x, u) = F_1(t,x, u) u . 
\end{equation} 
We look for solutions of \eqref{eq31} in lens shaped domains
\begin{equation}
\Omega_{r, \delta} = \big\{ (t, x) \ : \ t \ge  0 , \ \vert x - \ux \vert^2 + \delta t < r^2 \big\},   
\end{equation}
assuming that the equation is not hyperbolic at $(0, \ux, 0)$:
\begin{assumption}
There is $\uxi \in \RR^d $ such that the 
 matrix $\uA  := 
 \sum \uxi_j A_j(0, \ux, 0)$ has a nonreal eigenvalue.  
\end{assumption} 

The next theorem shows that the Cauchy problem is not 
H\"older well posed. We denote by $B_r$  the ball of radius $r$
centered at $\ux$. 

\begin{theorem}
\label{theo32}
For all $m$,  $\alpha \in ]0, 1]$, $r_0 > 0$  and $\delta > 0$, 
there are $r_\eps \to 0$, families of initial data $h_\eps \in H^m(B_{r_0})$ and solutions 
$u_\eps $ of $\eqref{eq31}$  on $\Omega_{r_\eps, \delta}$, such that 
\begin{equation}
\lim_{\eps \to 0}  \Vert u_\eps \Vert_{L^2(\Omega_{r_\eps, \delta})}  /
\Vert h_\eps \Vert^\alpha_{H^m(B_{r_0})}   = + \infty. 
\end{equation} 

\end{theorem}

\bigbreak

Let $\lambda_0$ denote an   eigenvalue of $\uxi \cdot A(0, \ux, 0)$ such that 
$\gamma_0 = \vert {\rm Im } \lambda_0 \vert > 0$ is maximum.  
Let $\ur $ denote an  eigenvector associated to $\lambda_0$. 
We consider  initial data
\begin{equation}
\label{eq35}
 h_{\eps}(x) := \eps^M \  {\rm Re } \big(  e^{i  x \cdot \uxi / \eps}   \ur \big) \,. 
\end{equation} 
We look for solutions
\begin{equation}
u_\eps(t, x) = \ru(t, x, t/\eps, x \cdot \uxi/\eps) 
\end{equation}
where $\ru (t, x, s, \theta)$ is $2 \pi$ periodic in $\theta$. For $u_\eps$ to be 
solution of the equation, it is sufficient that $\ru$ solves an equation of the form
\begin{equation}
\label{eq37}
\D_s \ru = \rA (y, \ru) \D_\theta \ru + 
\eps \big( \rB (y, \ru) \D_y u  + F(y, \ru) \big) ,
\end{equation}
with $y= (t,x-\ux)$ and 
$\rA (y, u) = \sum \uxi_j A_j (t, x, u)$. In particular
$ \rA (0, 0) = \uA$ and  the equation reads
\begin{equation}
(\D_s - \uA \D_\theta) \ru = \rF (\ru) := 
( \rA - \uA) \D_\theta \ru + 
\eps \big( \rB   \D_y u  + F(y, \ru) \big) . 
\label{eq37b}
\end{equation}
The solution of the Cauchy problem is given by 
\begin{equation}
\label{eq39}
\ru  = e^{ s \uA \D_\theta} \rh + \cT(\ru) \,, 
\quad \cT(\ru) (s) := \int_{0}^s  e^{ (s - {s'}) \uA \D_\theta}  \rF (\ru(s') ) ds'. 
\end{equation} 
We solve this equation  following the method    explained in Wagschall \cite{Wa}
(see also the references therein). 

\medbreak 
\noindent{\sl Function spaces and existence of solutions. }
Given power series 
$\ru = \sum u_\alpha y^\alpha$ and $\Phi = \sum \Phi_\alpha y^\alpha$, we say that 
$\ru \ll \Phi$ when $\vert u_\alpha \vert \le   \Phi_\alpha $ for all $\alpha$. 
Consider the series 
$$
\phi (z) =  c_0  \sum_{n=0}^{+\infty} \frac{z^n}{n^2 + 1} 
$$
where $c_0$ is taken such that $\phi^2 \ll \phi$ (cf \cite{La1}, \cite{Wa}). 
For $y \in \RR^{1+d}$, we denote by $Y = \sum {y_j}$, and we will consider 
power series $\ru(y)$ such that there is a constant $C$ such that
$$
\ru(y) \ll  C  \phi (R Y + R_0  ).   
$$
$R > 0 $ and $R_0 \in ]0, 1]$ are given parameters. These power series are convergent 
for $R \sum \vert y_j \vert + R_0 \le 1$.  

Next we introduce the weight function on $\ZZ$: 
$$
\langle n \rangle  = \vert n \vert \quad \mathrm{when} \ n \ne 0, 
\qquad \la 0 \ra = 2. 
$$
Note that for all $p$ and $q$ in $\ZZ$: 
\begin{equation}
\label{eq39}
\la p+q\ra \le \la p \ra + \la q \ra. 
\end{equation}

Given positive parameters $\gamma$, $\kappa$, $\eps$, $R$ and $\rho$, we consider 
formal Fourier series
\begin{equation}
\ru (s, \theta, y) = \sum_{- \infty} ^{+ \infty} \ru_n (s, y)  e^{i n \theta} 
\label{eq310}
\end{equation}
where the $\ru_n(s, y)$ are power series in $y$, with coefficient 
$C^\infty$  in $s \in [0,  \us  ]$, where 
\begin{equation}
\label{eq312}
\us :=  \min\{  \kappa/\gamma  , 1/ \eps \rho\}.
\end{equation}  
 We denote by  
  $\EE$  the space  of $\ru$   
such that there is a constant $C$ such that for all 
 $s \in [0,  \us  ]$: 
\begin{equation}
\ru_n( s, y) \ll C \frac{c_1}{n^2 + 1} e^{( \gamma s - \kappa)\la n \ra}
\phi (R Y + \eps \rho s) \,. 
\label{eq311}
\end{equation}
The number  $c_1$ is chosen such that 
$$
\sum_{p+ q = n} \frac{c_1}{p^2+1} \frac{c_1}{q^2+1} \le \frac{c_1}{n^2 + 1} . 
$$
Elements  $\ru \in \EE$ define smooth functions  on the domain 
\begin{equation}
\label{dom314}
\begin{aligned}
\Delta= \big\{ (s, \theta, y) \ : 
0 \le s < &  \us , \\
 & \theta \in \TT , \  R \sum \vert y_j \vert + \eps \rho s < 1 \big \}. 
\end{aligned}
\end{equation}
The best constant $C$ in \eqref{eq311} defines a norm $\VV u \VV$ on 
$\EE$. Equipped with this norm, $\EE$ is a Banach space. 
The choice of $c_0$ and $c_1$ and \eqref{eq39} imply that 
$\EE$ is a Banach algebra: 
\begin{equation}
\label{eq313}
\VV \ru \rv \VV \le \VV \ru \VV \ \VV \rv \VV . 
\end{equation}
 When $\ru$ is valued in $\CC^N$, we denote by 
 $\VV u \VV$ the sup  of the norms of the components of $\ru$.

 \begin{lemma}
 \label{lem32}
 If $F(y, u)$ is holomorphic on a neighborhood 
 of the origin in $\CC^{1+d} \times \CC^{N}$, there are constants
 $R_0$, $C_0$ and $a_0$ such that for all parameters 
   $\gamma$, $\kappa$, $\eps$, $R \ge R_0$ and $\rho$,
  the mapping $\ru \mapsto  F(\cdot, \ru)$  maps the ball of radius
  $a_0$ of $\EE$ into the ball of radius $C_0$ in $\EE$. 
 
 \end{lemma}
 
 \begin{proof}
 There are  constants  $R_0$, $C$  and $a$  such that
 $$
 F(y, u) \ll C \phi(R_0 Y) \prod_{j=1}^N \frac{1}{a - u_j}  ,
  $$
 in the sense of power series in $(y,u)$.
Substituting $u= u(y)$ in the expansion,  using \eqref{eq313} as well as the 
 identities $\phi^2 \ll \phi$ and $\phi (R_0 Y) \ll \phi(R Y + \eps \rho s)$, 
 yields for $\VV \ru \VV < a$: 
 $$
 \VV F( \cdot, \ru) \VV  \le C \frac{1}{(a - \VV \ru \VV)^N} . 
 $$
 \end{proof}

We further denote by $\VV \cdot \VV'$ the norm obtained 
when $\phi$ is replaced by its derivative $\phi'$ in \eqref{eq311}, 
and by $ \VV  \cdot \VV_1$ the norm obtained when 
$c_1/ (n^2+1)$ is replaced by $ c_1 / \sqrt{ n^2+1}$. In particular, there holds: 
\begin{eqnarray*}
\VV \D_y \ru \VV'  & \le  &  R  \VV \ru \VV,
 \\
\VV \D_\theta \ru \VV_1  & \le & \VV \ru \VV.  
\end{eqnarray*}
 Moreover, differentiating the estimate $\phi^2 \ll \phi $ implies that 
 $2 \phi \phi' \ll \phi'$, thus 
 \begin{equation*}
2  \VV \ru \rv \VV'  \le \VV \ru \VV \ \VV \rv \VV' .
 \end{equation*}
Similarly, there is $c_2$ independent of all the parameters such that: 
\begin{equation*}
  \VV \ru \rv \VV_1  \le  c_2 \VV \ru \VV \ \VV \rv \VV_1.
 \end{equation*}
Factoring out $y$ and $u$ in $\rA(y, u) - \uA$ and $u$ in $F(y, u)$, using that 
$y \ll ( 2 / c_0 R) \phi ( R Y ) \ll (2 c_0 / R) \phi (RY + \eps \rho s)$ and 
that $\phi' \ll \phi$, we deduce from 
Lemma \ref{lem32} and the estimates above that for 
$R \ge R_0$ and $\VV \ru \VV \le a_0$: 
\begin{eqnarray}
\label{eq316} \VV (\rA (y, \ru) - \uA) \D_\theta \ru \VV_1 & \le  &C
\big( R^{-1}  + \VV \ru \VV \big) \VV \ru \VV , 
\\
\label{eq317} \VV  \rB(y, \ru)  \D_y  \ru  + F(y, \ru) \VV'   & \le&  C
 R  \VV \ru \VV. 
\end{eqnarray}

Next, we investigate the action of the operator
$$
\rv (s) = \cI(\rf) (s)  := \int_0^s e^{ (s- {s'}) \uA \D_\theta }  \rf (s') ds' \,.
$$
On each Fourier component, it reads
$$
\rv_n(s, y ) = \int_0^s e^{  i  n (s- {s'}) \uA  }  \rf_n  (s', y ) ds'. 
$$
By the definition of $\gamma_0$, for all 
$\gamma > \gamma_0$ there is a constant $K_\gamma$ such that: 
\begin{equation}
\label{eq318}
\forall n \in \ZZ, \ \forall s \in [0, + \infty [: 
\quad e^{ i n s \uA } \le K_\gamma e^{ \vert n \vert  \gamma s } . 
\end{equation}
By definition of the norm $\VV \rf \VV_1$, there holds: 
$$
\rf_n (s', y) \ll \frac{c_1}{\sqrt{n^2 +1}} \VV \rf \VV_1 e^{ (s'  \gamma - \kappa) \la n \ra}
\phi( R Y + \eps \rho s') . 
$$
Using that for $s' \le s$, $\phi( R Y + \eps \rho s')  \ll \phi( R Y + \eps \rho s) $, and integrating 
term by term the power series in $y$, implies that 
$$
\rv_n (s, y) \ll \frac{c_1 K_{\gamma_1} }{\sqrt{n^2 +1}} \VV \rf \VV_1 
\phi( R Y + \eps \rho s) 
\int_0^s e^{ \vert n \vert \gamma_1 (s- {s'})} e^{ (s'  \gamma - \kappa) \la n \ra}
ds'
$$
For $\gamma > \gamma_1 > \gamma_0$, the last integral is estimated by 
$$
e^{ ( s \gamma - \kappa) \la n \ra} 
\int_0^s e^{ {s'} (\vert n \vert \gamma_1 - \la n \ra \gamma) } ds' 
\le  \frac {C}{ (\gamma- \gamma_1) \la n \ra }   e^{ ( s \gamma - \kappa) \la n \ra}  \,. 
$$
Choosing $\gamma_1 = (\gamma + \gamma_0) /2$, this shows that for 
all $\gamma > \gamma_0$, there is a constant $K_\gamma$ such that 
$$
\VV \cI( \rf) \VV \le K_\gamma \VV \rf \VV_1. 
$$ 
Similarly, there holds
$$
\rv_n(s, y) \ll \frac{K_ \gamma} {n^2 + 1} \VV \rf \VV' e^{ (s \gamma - \kappa) \la n \ra } 
\int_0^s e^{ \gamma (s - s') ( \vert n \vert - \la n \ra) }
\phi' (R Y + \eps \rho s') ds' . 
$$
Since $\vert n \vert \le \la n \ra$, we can ignore the exponential in the integral. 
Moreover, 
$$
\eps \rho  \int_0^s 
\phi' (R Y + \eps \rho s') ds'  \ll \phi (RY + \eps \rho s) - \phi(RY) \ll 
\phi(RY + \eps \rho s) . 
$$
Therefore, 
$$
\VV \cI( \rf) \VV \le \frac{K_\gamma}{ \eps \rho} \VV \rf \VV'. 
$$
Using \eqref{eq316} \eqref{eq317}, these inequalities yield estimates for the operator $\cT(\ru) $ defined in \eqref{eq39}. 
Similarly, one obtains estimates for increments $\cT(\ru)  - \cT(\rv)$:

\begin{prop}
There are $R_0$ and  $a_0$ and for all $\gamma > \gamma_0$ there 
is a constant $K_\gamma$, such that for all $R \ge R_0$, all $\kappa >0$, all
$\rho > 0$ and all $\eps \in ]0, 1]$, there holds for 
all $\ru$ and $\rv$  in $ \EE$  such that 
$\VV \ru \VV \le a_0$ and $\VV \rv \VV \le a_0$: 
$$
\begin{aligned}
&\VV \cT(\ru) \VV \le  K_\gamma \big( R^{-1} +  2 \VV \ru \VV + R \rho^{-1} \big) \VV \ru \VV, 
\\
&\VV \cT(\ru) - \cT(\rv) \VV \le  K_\gamma \big( R^{-1} + \VV \ru \VV
+ \VV \rv \VV   + R \rho^{-1} \big) 
\VV \ru - \rv  \VV
\end{aligned}
$$

\end{prop}

\begin{cor}
\label{cor35}
With notations as above, if 
$$
K_\gamma \big( R^{-1} +  4 a + R \rho^{-1} \big)  < \frac{1}{2}, 
$$
then for all  $\rf \in \EE$ with $\VV \rf \VV \le a$, the equation 
$$
\ru = \rf + \cT(\ru)
$$ 
has a unique solution $\ru \in \EE$ such that 
$\VV \ru \VV \le 2 a$. Moreover, 
$$
\VV \ru - \rf \VV \le K_\gamma \big( R^{-1} +  4 \VV \rf \VV  + R \rho^{-1} \big)
\VV \rf \VV. 
$$

\end{cor}

\bigbreak

\noindent{\sl Application.} In accordance with \eqref{eq35}, we solve the Cauchy problem
\eqref{eq37} with initial data
\begin{equation}
\label{eq319}
\ru_{\vert s = 0} = \rh := \eps^M \mathrm{Re}  ( e^{ i \theta} \ur ) . 
\end{equation}
Let 
\begin{equation}
\rf = e^{s \uA \D_\theta} \rh  = \eps^M \mathrm{Re}  \big( e^{i s \lambda_0 +  i \theta} \ur \big)  . 
\label{eq320}
\end{equation}
We consider only small values of the  parameter $\eps$, and we use the notation 
$\eps^M = e^{- \kappa_1}$, that is $\kappa_1 =  M \vert \ln \eps \vert$.

Consider a small parameter $\beta > 0 $, to be chosen later on,  such that 
$\beta M <  1/2 $. We fix 
\begin{equation}
\label{eq321}
\left\{ \begin{array}{ll}
  \gamma = (1 + \beta) \gamma_0, \quad  & \kappa = (1- \beta) \kappa_1, \quad
\\
  R = e^{ \beta \kappa_1} = \eps^{ - \beta M} , \quad & \rho = R^2 = \eps^{- 2 \beta M}. 
\end{array}\right.
\end{equation}
Introduce $\sigma = (1 - \beta)/ (1 + \beta) < 1$. 
For $\eps$ small enough, $\kappa/ \gamma = \sigma \kappa_1/\gamma_0 \le 
(\eps \rho)^{-1} = \eps^{ -1 + 2 \beta M}$, thus, the end point \eqref{eq312}
is $\us = \sigma \kappa_1/ \gamma_0 $. 

\begin{prop}
\label{prop36}
There is a constant $c > 0$ such that 
for all $M\ge 1 $ and $\beta \in ]0, 1/ 2M[ $, there  is $\eps_0$ such that for all $\eps \in ]0, \eps_0]$,
and parameters as in $\eqref{eq321}$, the Cauchy problem 
$\eqref{eq37}, \eqref{eq319}$ has a solution 
$\ru \in \EE$, and 
\begin{equation}
\label{est322}
\forall (s, \theta, y) \in \Delta \, : \quad
\vert \ru(s, \theta, y) \vert \ge   c e^{ s \gamma_0 - \kappa_1} . 
\end{equation}
\end{prop}

\begin{proof}

For $\eps$ small enough, there holds
$$
\VV \rf \VV = \frac{2}{c_1 c_0}  \vert \ur \vert  \max_{s \in [0,  \us[} 
e^{ \kappa - \kappa_1 + s (\gamma_0- \gamma) } \le C e^{ - \beta \kappa_1}
$$
for some constant $C$ independent of $\eps$ and $\beta$. 
By Corollary \ref{cor35}, there is $K$, depending only on $\beta$, such that
for $ K    \eps^{ \beta M} < 1$, the problem has a unique solution
$\ru \in \EE$ and 
$$
\VV \ru - \rf \VV \le K  e^{ - 2 \beta \kappa_1}. 
$$
 For $(s, \theta, y) \in \Delta$, there holds
$R \sum \vert y_j \vert + \eps \rho s \le 1$. Since the series 
$\phi (z)$ converges at $z =1$, 
 $$
 \vert (\ru - \rf) (s, \theta, y) \vert 
 \le  K e^{ - 2 \beta \kappa_1} \sum_{n \in \ZZ} \frac{c_1}{n^2 +1} \phi(1) 
 e^{ ( s \gamma - \kappa)\la n \ra} .
 $$
 Since $s \gamma - \kappa \le 0$ and $\la n \ra \ge 1$, this implies that there is 
 $K'$ such that 
$$
\begin{aligned}
 \vert (\ru - \rf) (s, \theta, y) \vert 
   &  \le  K'   
e^{ - 2 \beta \kappa_1}  
 e^{ ( s \gamma - \kappa) }
 \\  
  &  \le  K'  e^{ s \gamma_0- \kappa_1} 
e^{ -  \beta \kappa_1}    e^{  \us \beta \gamma_0}  = 
 K'  e^{ s \gamma_0- \kappa_1} 
e^{ -  \beta (1- \sigma) \kappa_1 }    . 
\end{aligned}
 $$
Because $\ur$ and $\overline {\ur}$ are eigenvectors associated to distinct eigenvalues 
$\lambda_0$ and $\overline{ \lambda_0}$, they are linearly independent and there is 
$c > 0$ such that 
$$
\vert f( s, \theta) \vert   \ge 2 c e^{ s \gamma_0 - \kappa_1} . 
$$
Since $\sigma < 1$, the two estimates above imply that for $\eps$ small enough \eqref{est322}
is satisfied. 
\end{proof}

\bigbreak
\begin{proof}[Proof of Theorem $\ref{theo32} $]

The integer $m\ge 1$ and the H\"older exponent $\alpha\in ]0, 1] $ are given, as well as 
the parameter $\delta > 0$. 
We fix $M$ large enough, such that 
\begin{equation}
\label{rel323}
\alpha' := \frac{M - m}{M} \alpha - \frac{ 1 + d}{ 2 M } > 0  . 
\end{equation} 
Note that $\alpha' < \alpha \le 1$. Next we choose $\beta > 0$ such that 
\begin{equation} 
\label{rel324}
1 - \alpha' < \sigma := \frac{1- \beta}{1+ \beta} \quad \mathrm{and} \quad
2 M \beta < 1  
\end{equation}
and we fix the parameters $\gamma$, $\kappa$, $R$ and $\rho$ as in 
\eqref{eq321}. 
By Proposition \ref{prop36}, for $\eps$ small enough, we have a solution 
$\ru$ of \eqref{eq37} \eqref{eq319} on the domain $\Delta$ defined in \eqref{dom314}. 
Thus 
$$
u_\eps (t, x) = \ru(\frac{t}{\eps}, \frac{x \cdot \uxi}{\eps}, t, y)
$$ 
is a solution of \eqref{eq31} \eqref{eq35} on the domain 
\begin{equation*}
\tilde \Delta_\eps = \big\{ (t, x) \ : 
0 \le t  <     \eps \us ,  \  
    \   \sum \vert x_j \vert +   t   +   t \eps^{ - \beta M}      <  \eps^{ \beta M}  \big \}. 
\end{equation*}
Since $\ut_\eps := \eps \us = \sigma \gamma_0^{-1} M  \eps \vert \ln \eps \vert $ and 
$2 \beta M < 1$,  for   $\eps$ small enough,  this domain contains 
\begin{equation*}
  \Delta_\eps = \big\{ (t, x) \ : 
0 \le t  <     \ut_ \eps  ,  \  
    \       \vert x  \vert       <   c \eps^{ \beta M}  \big \} 
\end{equation*}
with $c = 1/(2 \sqrt d)$. For $\eps $ small, it also contains the lens shaped domain 
$\Omega_{r_\eps, \delta} $, with 
\begin{equation}
\label{def325}
r_\eps = (t_ \eps/ \delta)^{ 1/2} . 
\end{equation}

Moreover, for $\eps$ small, $\Omega_{r_\eps, \delta} $ contains
the cube 
$$
\big\{ \ut_\eps - \eps  \le t \le \ut_\eps, \ \vert x - \ux \vert \le \eps \big\} \,. 
$$
Thus, Proposition \ref{prop36} implies that there is $c > 0$ such that for all 
$\eps $ small enough: 
$$
\Vert u_ \eps \Vert_{L^2(\Omega_{r_\eps, \delta})} \ge c e^{ \us \gamma_0 - \kappa_1}
\eps^{1+ d}  =  c \eps^{ M ( 1 - \sigma) + (1 + d)/2}. 
$$

On the other hand,  the Sobolev norm of the initial data 
on a fixed ball $B_{r_0}$ centered at $\ux$ is of order: 
\begin{equation*}
\Vert h_{\kappa, \eps } \Vert_{H^m (B_{r_0})}  \le C   \eps^{ M - m} \,. 
\end{equation*}
Thus, using the notation \eqref{rel323}, 
\begin{equation}
   \Vert u_\eps \Vert_{L^2(\Omega_{r_\eps, \delta})}  /
\Vert h_\eps \Vert^\alpha_{H^m(B_{r_0})}    \ge \frac{c}{C^\alpha} 
\eps^{ M (1 - \sigma - \alpha')} 
\end{equation} 
which, by \eqref{rel324}, tends to  $+ \infty$ as $\eps$ tends to zero. 
\end{proof}


\section{Solvability in Sobolev spaces}

In this section, we consider the fully nonlinear Cauchy problem  in $\RR^{1+d}$: 
\begin{equation}
\label{eq41}
\left\{\begin{aligned} 
 & \D_t u = F(t, x, u, \D_{x_1} u, \ldots, \D_{x_d} u ) , \quad  t \ge 0
\\
& u_{\vert t = 0} = h,
\end{aligned}\right. 
\end{equation}
near $(0, \ux)$. We assume that $F$ is $C^\infty$ in a neighborhood of $\up := (0, \ux, \uu, \uv)$ 
in $\RR\times \RR^d \times \RR^N \times (\RR^N)^d$. The initial data $h$ is smooth
and satisfies
\begin{equation}
\label{eq42}
 h(\ux) = \uu , \quad   \D_x h(\ux)= \uv.  
\end{equation}

\begin{assumption}  
\label{ass41}
 There is $\uxi \in \RR^d$ such that the  matrix 
  $\uxi \cdot \partial_v F ( \up )= \sum \uxi_j \D_{v_j} F(\up) $ has   nonreal eigenvalues. 
\end{assumption}

Because $\xi \cdot \D_v F$ is real,  this implies that there is at least one eigenvalue 
   with positive imaginary part.  Denote by  $\underline \Pi $  the spectral 
  projector of 
$\uxi \cdot \partial_v F ( \up )$ 
associated to   eigenvalues  in $ \{ \mathrm{Im} \mu > 0 \}$.

\begin{theorem}
 \label{theo42}
Let     $s  > d/2+1$ and   $ s \le s'  <  2 s - 1 - d/2 $. 
Suppose that the Cauchy data $h$ satisfies $\eqref{eq42}$ and 
\begin{equation}
\label{b43}
(\Id - \underline \Pi) h \in H^{s'} \quad near \ \ux. 
\end{equation} 
If  
the Cauchy problem $\eqref{eq41}$ has a solution in 
$C^0([0, T]; H^{ s }(\omega))$, then  
\begin{equation}
\label{b44}
 (\ux, \uxi) \notin WF_{H^{s'}} ( \underline \Pi  h). 
\end{equation} 
\end{theorem}

In this statement, $WF_{H^{s'}}$ denotes the $H^{s'}$ wave front set. 
In the spirit of Hadamard's argument and of Theorem \ref{th21}, Theorem \ref{theo42} shows 
that smoothness of part of the Cauchy data, here $\underline \Pi h$, implies 
smoothness of the other components.  
For all $s'' \in [s, s' [$, there are many Cauchy data $h$ such that 
\begin{equation}
\label{b45}
 h \in H^{s''}, \quad (\Id - \underline \Pi) h \in H^{s'} , \quad  
  (\ux, \uxi) \in WF_{H^{s'}} ( \underline \Pi  h). 
\end{equation} 
For such data, all $T >0$ and all neighborhood $\omega$ of $\ux$, 
Theorem \ref{theo42} implies  that 
 \eqref{eq41} has no local solution 
in $C^0([0, T]; H^{ s }(\omega))$. 
This implies that the Cauchy problem is not locally well posed from
$H^{s''}$ to $C^0(H^s)$ for all 
$s'' = 2 s -1 - d/2 - \eps > s$ when $s > 1 + d/2$.

 The proof is an application of the results of 
Monique Sabl\'e-Tougeron \cite{ST}  about the propagation of microlocal singularities 
for nonlinear boundary value problems. 
For the convenience of the reader, we sketch a proof 
within the class of spaces  $C^0(H^s)$  instead of  the class $H^{s,s'}(\RR^{1+d})$
used in \cite{ST}.

\begin{proof}
   Decreasing slightly 
$s$, we can assume that  
$\rho := s -1 - d/2 \notin \NN$. Suppose that  $h \in H^{s}(\RR^d)$ satisfies
\eqref{eq42} and that 
 $u \in C^0([0, T];  H^{s}(\omega))$  is a solution of \eqref{eq41}.

{\bf a)} 
The product is   continuous from 
$H^{\sigma- \alpha} \times H^{\sigma- \beta}$ into  $H^{\sigma - \alpha - \beta}$, 
when   $\sigma > d/2$, $\alpha \ge 0$, $\beta \ge 0$ and
$\alpha + \beta \le 2 \sigma$.   
By induction on $k$, \eqref{eq41} implies that  
\begin{equation*}
\partial_t^k u \, \in \, C^0 ( [0; T] H^{ s - k }(\omega))   \,,   \quad k \in \{0, \ldots, [2 s - 2] \}\, 
\end{equation*}
Therefore, for all smooth function  
$G  $,  
\begin{equation*}
\partial_t^k G(\, \cdot \, , \,\cdot \, , u, \partial_x u)  \, \in \, C^0 ( H^{ s- 1 - k })  \, ,  \quad k \in \{0, \ldots,
[2 s - 2]
\}\, .
\end{equation*}
Since,  $\rho = 2s - 2  -  (s - 1  + d/2) < 2s -2 - d$ this 
property is true up to $ k = [\rho ] + 1$. 
Denoting by   $C^\alpha(\RR^d)$ the usual H\"older space 
for  $\alpha \in \RR \backslash \ZZ$,  
this implies that $ g  =  G( \cdot, u(\cdot), \D_x u (\cdot))$ satisfies 
\begin{equation} 
\label{eq43}
\partial_t^k g   \, \in \, C^0 ([0, T];  C^{ \rho - k }(\omega))  \, ,  
\quad k \in \{0, \ldots, [\rho]+ 1 \}\, .
\end{equation}
For  $\rho > 0 $, $\rho \notin \NN$, we denote by  $\widetilde C _\rho$ the set of functions $g$ which satisfy this 
 property. 
 \medbreak

{\bf b)} Localizing near $(0, \ux)$, and using Bony's paralinearization theorem
 in $x$ (\cite{Bo}, \cite{ST}),  
 \eqref{eq41} implies that $\tilde u = \chi_1 u$ satisfies 
 \begin{equation}
 \label{eq44}
\partial_t \tilde u \, -  \, T_A (t, x, \partial_x) \, \tilde u \, = \, f ,  
\end{equation}
where $\chi_1 \in C_0^\infty(\omega)$ is equal to one 
 near $\ux$ and   $f \in C^0([0, T]; H^{ s - 1 + \rho}(\omega')) $ 
 for some smaller neighborhood $\omega'$ of $\ux$.  
 In this equation, $T_A$ denotes a  paradifferential operator 
 in $x$ of  symbol 
\begin{equation} 
A (t, x, \xi):= 
\sum_{j=1}^d \xi  \cdot  \partial_v F (p(t, x) )  +  \partial_u F p(t, x)
\end{equation}  
 with  $ p(t, x) =  (t, x, u(t, x),  \partial_x u(t,x))$. The coefficients  belong to
  $\widetilde C^\rho$.

\medbreak 

{\bf c)}  We perform a microlocal block diagonalization of $A$. 
 Near $(0, \ux, \uxi)$,  there are symbols 
\begin{equation*}
P \, = \, \sum_{j = 0} ^{ [\rho]} P_{j }\, , \quad
D \, = \, \sum_{j = 0}^{ [\rho]} D_j\,
\end{equation*}
  such that  
\begin{equation}
\label{eq46}
P \sharp A  -  \partial_t (P - P_{[\rho]}) \,  = \,  D \sharp P\,. 
\end{equation}
The terms  $P_j(t, x, \xi)$ 
[resp.  $D_j$]  are  $C^\infty$ and homogeneous of degree $   -j$
[resp. $1-j$] in  $\xi$ and $\widetilde C^{\rho - j}$ in $(t,x)$. 
 In \eqref{eq46},  $\sharp$ denotes the composition of symbols (cf \cite{Bo}): 
 $$
\Big(  \sum_{j \le [\rho]} P_j \Big) \sharp 
\Big( \sum_{j \le [\rho]} Q_j\Big) \, := \, 
\sum_{ l + m +\vert \alpha \vert \le [\rho]} 
 \  {1\over \alpha!} (\partial_\xi^\alpha P_l) \, ( (- i \D_x)^\alpha Q_m)
$$
Moreover, $D$ is block diagonal  and there is one block $D^I$  associated to 
  the spectrum of $ A(0, \ux, \uxi)$ in 
$\{ \mathrm{Im} \mu > 0 \}$. In particular 
\begin{equation}
\label{eq47}
\mathrm{Im \ spec} (D^I_0(0, \ux, \uxi)  )   > 0
\end{equation}  
The construction is classical. The principal terms $P_0$ and $D_0$ are chosen such that 
$$
P_0 A_0 (P_0)^{-1} = D_0,
$$
where $A_0 = \xi \cdot \D_v F$ is the principal symbol of $A$. 
Next, one proceeds by induction on $j$, choosing $P_j$ and $D_j$ such that 
the terms of degree $1-j$ in the two sides of   \eqref{eq46}  are equal.  
In particular,  \eqref{eq46} is an identity between symbols of degree $1$, 
like $D$,  and the degree of the last term, like $D_{[\rho]}$ is 
$1 - [\rho]$. This is why, in  $\partial_t P$,  we can ignore the last term which would be 
$\partial_t P_{[\rho]}$, of degree  $- [\rho]$.   

In this computation, we only use symbols with positive degrees of smoothness  
$\rho -j$, with $j \le [\rho]$.  
  However, the term  $\partial_t P_{[\rho]}$ 
which will appear in the remainders, requires one more derivative. 
This is why we  took  $k \le [\rho]+1$ in \eqref{eq43}.
 
\medbreak 
{\bf d)}   Suppose that  $\chi  (t, x , \xi) $  is a microlocal cut-off 
function supported  in a conical neighborhood of $(0, \ux, \uxi)$ where 
 \eqref{eq46} is satisfied. 
The equation \eqref{eq44} implies that 
\begin{equation*}
(\partial_t  -  T_ {\chi  D} )   T_{\chi P} \tilde u \, = \, 
T_{ \chi^2 P} (\partial_t - T_{A} \tilde u)  \, +  \, T_{ \chi \partial_t P_{[\rho]}} \tilde u \, + \,  
T_{ Q} \tilde u \,  +  \, R u
\end{equation*}
where  $R$  is a remainder in the $[\rho]$-calculus   $x$. If 
$\rho < 1$ then  $Q = 0$ and if  $\rho > 1$, $Q$  is a symbol of degree zero, 
  $(\rho-1)$ smooth  in $x$ and equal to zero near  $(0, \ux , \uxi)$.  
In particular, near $(0, \ux)$, 
 $R u \in C^0( H^{s -1 + \rho})$  
and  $T_ Q u (t, \cdot) \in  H^{s -1 + \rho})$ near  $(\ux, \uxi))$, uniformly in  $t$.  
Moreover, since  $\partial_t P_{[\rho]}$ is of degree $- [\rho]$ 
 with smoothness $C^{\rho - [\rho]-1}$
in  $x$, the operator  $ T_{ \chi \partial_t P_{[\rho]}}$ is of order 
$- [\rho]+ (1 -  \rho +[\rho]) = 1 - \rho$.  
 Therefore, we  see that 
 $$
w := T_{\chi P} \tilde u 
$$
satisfies 
\begin{equation}
 \label{eq48}
\partial_t  w  - T_{\chi D}    w  \in C^0( H^{ s-1 + \rho} ) \,, \quad
{\rm near   }\   (0, \ux, \uxi)\, . 
\end{equation} 
\medbreak 

{\bf e)}   $D$ is block diagonal. Denote by   $w_I$ the components of 
$w$ which correspond to the bloc $D^I$.  The equation \eqref{eq48} implies that 
 \begin{equation}
 \label{eq49}
\partial_t  w_I  - T_{\chi D^I}    w  \in C^0( H^{ s-1 + \rho} ) \,, \quad
{\rm near   }\   (0, \ux, \uxi)\, . 
\end{equation} 
 By \eqref{eq47}, this problem is elliptic and the backward Cauchy problem is well posed
 (\cite{ST}).  This implies that  $ w  \in C^0( H^{ s + \rho} )$ 
near  $(0, \ux, \uxi)$. 
 By construction,  
$w_I = T_\Pi u $  where   $ \Pi = \sum_{j \le [\rho]} \Pi_j $, with $\Pi_j$ of 
  degree $-j$ in $\xi$ and $\widetilde C^{\rho -j}$ in $(t,x)$.    
  In particular, 
\begin{equation}
  \label{eq49b}
w_I {}_{\vert t = 0} =  T_{\Pi_{\vert t = 0} }  h   \in   H^{ s + \rho}   \, , 
\, \quad
{\rm  near  }\   (\ux, \uxi)\, . 
\end{equation}
 For  $(t,x, \xi)$ close to 
$(0, \ux, \uxi)$, the  principal symbol $\Pi_0 (t, x, \xi) $ is the spectral 
  projector of 
$\xi \cdot \partial_v F (t,x, u(t,x), \partial_xu(t,x) )$ 
corresponding to eigenvalues in $\{ \mathrm{ Im } \mu > 0 \}$. 
In particular  
  \begin{equation}
  \label{eq410}
   \underline  \Pi  := \Pi_0(0, \ux, \uxi). 
    \end{equation}
    Since the system $(\Pi_0, \Id - \underline \Pi ) $ is elliptic near $(\ux, \uxi)$, 
    there are symbols $U = \sum_{j \le [\rho]} U_j$ and 
    $V_j = \sum_{j \le [\rho]}$ of degree zero such that 
 $$
\underline \Pi  \, = \, U \, \sharp \, \Pi _{\vert t = 0} \, + \, V \, \sharp \, (\Id  - \underline \Pi), 
 \quad \mathrm{near} \ (\ux, \uxi) \, . 
$$
This implies that  for $\chi_1$ supported in a sufficiently small conical neighborhood
 of $(\ux, \uxi)$
\begin{equation}
 \label{eq411}
\chi_1(x, D_x) \underline \Pi   =  T_U   T_{\Pi_{\vert t = 0}}   +   T_V   (\Id - \underline \Pi)  +   R 
\end{equation}
with  $R$ of order  $-\rho$.   
 
 Suppose that the initial data satisfies 
\begin{equation} 
 \label{eq413} 
h \in H^{s}(\RR^d)  \, , \quad (\Id - \underline \Pi ) h \in H^{s + \rho} (\RR^d) . 
 \end{equation}
 Then 
 \eqref{eq411} and \eqref{eq49b}   imply that 
 $ \chi_1 (x, D_x) \underline \Pi h  \in H^{s+\rho}$, 
that is 
$ \underline \Pi  h \in  H^{s+ \rho}$ near  $ (\ux, \uxi)$ or 
$(\ux , \uxi) \notin WF_{H^{s+ \rho}} (\underline \Pi h)$. 
\end{proof}
 
\begin{remark} 
a) Note that only the condition $u \in C^0( H^s) $ is used to prove \eqref{eq49b}.  
  
  b)  The proof only relies on the ellipticity 
$\partial_t - T_{\chi D_I}$. Thus, the Sobolev spaces 
$H^s$ do not play any particular role and there are analogous results  
in the   H\"older spaces $C^\mu$.

  c) For semilinear equation, the critical index  $1 + d/2$ 
can be decreased to $d/2$ as usual. 
This is also the case for  system of conservation laws,  since we only need
to paralinearize functions of $u$. 
 \end{remark}

\bigbreak 

One can push a little further the analysis  when Assumption \ref{ass41} is strengthened.

\begin{assumption}   
\label{ass51}
The real eigenvalues of 
  $\xi \cdot \partial_v F ( p(t,x))$  are semi-simple and have constant multiplicity, 
and there are nonreal eigenvalues.  
\end{assumption}

In this case, the condition $s' < 2s - 1 - d/2$ in Theorem \ref{theo42} can be relaxed.

 \begin{theorem} Under Assumption $\ref{ass51}$,
 for all   $\sigma > d/2+1$, there are Cauchy data 
  $h \in H^{\sigma } (\RR^d)$, satisfying $\eqref{eq42} $ such that 
 for all  $s  > 2+ d/2$, all $T> 0$ and all neighborhood
$ \omega$ of $\ux$, the Cauchy problem $\eqref{eq41}$ has 
 no solution 
 $u \in C^0([0, T]: H^s(\omega))$. 
  \end{theorem}

\bigbreak 
The meaning is that one can take $\sigma$ very large and 
  $s$ very close to $ 2+d/2+ 2$,  so that $u$ 
  will be of class $C^2$, but not much smoother, while the initial data is 
  as smooth as we want.  
 
\bigbreak 
\begin{proof}
{\bf a)} Suppose that  $u \in C^0([0, T]; H^s(\omega))$ solves \eqref{eq41}.  
  We show that 
 $u \in C^0([0, T'];  H^{s'}(\omega'))$  for $T' < T$, 
$  \omega' \subset\!\subset \omega$ and  $s'\le \sigma $ such that 
$s' <  2s -2 - d/2  $. 
It is sufficient to prove that  $u \in C^0([0, T'];  H^{s'}(\omega'))$
with $s' =\min (\sigma,  2s - 2 - d/2)$ when 
$\rho := s- 1 - d/2 \notin \NN$. 

The analogue  is proved in \cite{ST}, for $m$-th order scalar equations, when 
the real roots of the principal symbol are simple. 
As in the proof of Theorem \ref{theo42}, near any $\xi \ne 0$, there is an elliptic symbol 
$P = \sum _{j \le [\rho]}$ such that 
$ w := T_{  P} u  $ 
satisfies 
$$
\partial_t  w  - T_{  D}    w  \in C^0( H^{ s'}  ) \,, \quad 
w_{\vert t = t_0} \in H^{\sigma} \, \quad
{\rm   near   }\   (0, \ux, \xi )\, .  
$$
 The matrix $D$ is block diagonal. By Assumption \ref{ass51}, the blocks of 
the principal symbol $D_0$ are either hyperbolic, that is of the form $ i \lambda \Id$ with 
$\lambda(t, x, \xi)$ real, or elliptic, meaning that the imaginary part of the eigenvalues 
is either positive or negative. 

Since the Cauchy data is $H^\sigma$, the equation implies that  hyperbolic blocks
are microlocally $H^{s'}$. The same result holds for negative elliptic blocks. 
For positive elliptic block, we use the backward elliptic
regularity,  and decreasing the interval of time, 
we see that the elliptic modes are $C^0(H^{s + \rho})$. This shows that  
$w\in C^0(H^{s '}) $.  Since $P$ is elliptic,   $u$ has the  
same regularity.

\medbreak

{\bf b)} Repeating the argument in  a), we deduce that any solution  
$u\in C^0(H^s)$ with initial data in $H^\sigma$  is necessarily in  $C^0(H^\sigma)$ on a smaller domain. 
Therefore, by Theorem \ref{theo42}, if 
 $$ 
h \in H^{\sigma}(\RR^d)  \, , \quad (\Id - \underline \Pi ) h \in H^{\sigma' + \rho'} (\RR^d) \,, \quad
 \underline \Pi_0 h \notin H^{\sigma'+ \rho'} (\ux, \uxi)\, ,
$$ 
the Cauchy problem has no solution in 
$C^0(H^{\sigma'})$, thus no solution  $u\in C^0(H^s)$. 
  \end{proof}
\bigbreak


\section{An example of problems with elliptic  zones}

   In this section we   consider    a modified version of \eqref{ex11}. 
   This is a nonhyperbolic form of Kirchhoff equation.  
   The advantage is that  we can make explicit computations, the drawback is
   that the equation is nonlocal\footnote{Thierry Colin 
   proposed the   simpler   example
 :  $\partial_t u \, =  - (1 - \Vert \partial_x u \Vert^2_{L^2}) \, \partial_x^2 u $. 
The system \eqref{eq56} is  first order and fits  the general 
 presentation of this paper. }.
  The modified system reads: 
\begin{equation}
\label{eq56}
\left\{\begin{aligned}
&\partial_t u  =  a(t)  \partial_x v  , 
\cr
& \partial_t v = \vert a(t) \vert  \partial_x u  ,  
\end{aligned}\right. \qquad  \mathrm{with}\quad 
a(t) = \Vert u(t) \Vert^2_{L^2}  -1 . 
\end{equation}
 As \eqref{ex11}, this system has a natural (formal) energy:   
 \begin{equation}
 E (t) = \Vert v(t, \cdot)  \Vert^2_{L^2} + \Big\vert  \Vert u(t, \cdot)  \Vert^2_{L^2} - 1  \Big \vert . 
 \end{equation}
If $u$ is $C^1 (L^2)$, the mapping $t \mapsto U(t) :=  \Vert u(t, \cdot)  \Vert^2_{L^2}$ is 
$C^1$, thus  $\vert U(t) - 1 \vert$ is  Lipschitzian and 
$$
\frac{d}{dt}  \vert U(t) - 1 \vert = \mathrm{sign} \big(U(t)- 1\big) \frac{d U(t)}{dt}  
\quad  a.e. 
$$
If $(u, v)$ is in addition $C^0(H^1)$, then 
$$
\frac{d E(t) }{dt}   = \vert a \vert  \int v \D_x u dx + a \mathrm{sign}(a) 
\int u \D_x v  dx = \vert a \vert \int \D_x (u v) dx = 0
$$

In the spirit of \eqref{ex14}, we considered the filtered system, 
with truncated frequencies. Since the system has constant coefficients in $x$, 
it is sufficient to filter the initial data:  
\begin{equation}
\label{eq57}
\left\{\begin{aligned}
&\partial_t u^\lambda   =  a^\lambda (t)  \partial_x v^\lambda  , 
\cr
& \partial_t v^\lambda  = \vert a^\lambda (t) \vert  \partial_x u^\lambda   ,  
\end{aligned}\right. \qquad  
\left\{\begin{aligned}
& u^\lambda{}_{\vert t = 0}  = 0  , 
\cr
&   v^\lambda {}_{\vert t = 0}   =  S_\lambda h   ,  
\end{aligned}\right. 
\end{equation}
with $a^\lambda (t) = \Vert u^\lambda (t) \Vert^2_{L^2}  -1$ and $S_\lambda$ 
is defined on the Fourier side by:  
   $$
\widehat{S_\lambda  h } (\xi)  =  1_{\{ \vert \xi \vert \le \lambda\} } \hat h(\xi). 
$$
For Fourier transforms, the system reads  for 
$\vert \xi \vert \le \lambda$:
\begin{equation}
\label{eq58}
\left\{\begin{aligned}
&\partial_t \hat u^\lambda   =  i \xi  a^\lambda (t)   \hat v^\lambda  , 
\cr
& \partial_t \hat v^\lambda  =  i \xi  \vert a^\lambda (t) \vert   \hat u^\lambda   ,  
\end{aligned}\right. \qquad  
\left\{\begin{aligned}
& \hat u^\lambda{}_{\vert t = 0}  = 0  , 
\cr
&   \hat v^\lambda {}_{\vert t = 0}   =  \hat  h   ,  
\end{aligned}\right. 
\end{equation}
and $\hat u^\lambda = \hat v^\lambda = 0$ for $\vert \xi \vert \ge \lambda$. 
This is a system of ordinary differential equations, and it has local  
solutions, $C^1$ in time with values in $L^2$. One can use the energy $E$
to prove that the solutions are global in time, but we  provide a direct proof. 

Suppose that $(u^\lambda, v^\lambda)$ is defined 
and that $U^\lambda(t) := \Vert u^\lambda(t) \Vert^2_{L^2} \le 1$  on 
$[0, T]$. This is certainly true for $T$ small. Then,  
$\vert a^\lambda \vert = - a^\lambda$ on this interval and for $\vert \xi \vert \le \lambda$: 
\begin{equation}
\label{eq59}
\left\{\begin{aligned}
&  u^\lambda (t, \xi )   =   i \sinh ( \xi A^\lambda(t) \big) h(\xi) 
\\
&   v^\lambda (t, \xi)  = \cosh ( \xi A^\lambda(t) \big) h(\xi)  
\end{aligned}\right. \qquad  
 A^\lambda (t) = t -  \int_0^t U^\lambda(s) ds. 
\end{equation}
Therefore 
$$
U^\lambda(t)  = \frac{1}{2\pi} \int_{\vert \xi\vert \le \lambda} 
  \sinh^2 (\xi  A^\lambda(t))   \,
 \vert \hat h (\xi)\vert^2 \, d\xi , 
$$
and 
\begin{equation}
 \frac{ d U^\lambda }{dt }  = (1 - U^\lambda(t))  I^\lambda  (t)  \, 
\label{eq510}
\end{equation}
with 
$$
I^\lambda (t) = \frac{1}{2\pi} \int_{\vert \xi\vert \le \lambda} 
 \xi  \sinh ( 2 \xi  A^\lambda(t))   \,
 \vert \hat h (\xi)\vert^2 \, d\xi . 
$$
Since $U^\lambda(0) = 0$, the equation \eqref{eq510} implies that $1 - U^\lambda$ does not 
vanish and remains positive on $[0, T]$.  In particular 
$U^\lambda (T) < 1$. By \eqref{eq59}, there holds
$$
\Vert v^\lambda (t) \Vert^2_{L^2} = U^\lambda(t) + \Vert S_\lambda h \Vert^2_{L^2}
\le 1 + \Vert h \Vert^2_{L^2}. 
$$
Therefore, by continuation, this implies that \eqref{eq57} has a unique global solution
in $C^0([0, + \infty[ ; L^2(\RR)$ and that 
$U^\lambda(t) < 1$ for all time. 
Moreover,   $A^\lambda \ge 0$   and the integral 
$I^\lambda $ is positive, implying that $U^\lambda$ is strictly increasing.

 Since $U^\lambda$ is increasing, there holds
 $A(t) \le t (1 - U^\lambda(t))$. Therefore 
$$
 \begin{aligned}
1 \ge U^\lambda(t) &  \ge  
 \frac{1}{8  \pi} \int_{C_ \lambda } 
\big(  e^{2 \vert \xi \vert t(1 - a^\lambda(t))  }  - 2 \Big) \,
 \vert h (\xi)\vert^2 \, d\xi  
 \\ 
 & \ge  
e^{ t(1 - U^\lambda(t)) \lambda }  \int_{C_\lambda }    \frac{1}{8  \pi}   
 \vert h (\xi)\vert^2 \, d\xi   -  \frac{1}{2} \Vert h \Vert^2_{L^2}. 
\end{aligned}
 $$
where  $C_\lambda := \{ \lambda  2  \le \vert  \xi\vert \le \lambda\}$. 
 Hence 
\begin{equation}
 \label{eq510b}
1 - U^\lambda(t) \, \le \, \frac{1} {  t \lambda  } ( \mu(\lambda ) + K  ) \, .
\end{equation} 
 where   $K  = \ln ( 8 \pi (1 + \Vert h \Vert^2)) $ and 
 \begin{equation}
 \mu( \lambda) = - \ln \Big(  \int_{C_\lambda }   \vert h (\xi)\vert^2 \, d\xi  \Big) 
 \end{equation}
 The condition $\mu (\lambda) \le C \lambda$ implies that 
 $h$ is real analytic. On the other hand, for general non analytic functions, 
 there holds 
 \begin{equation}
 \label{eq511}
 \lim_{ \lambda \to + \infty} \frac{ \mu (\lambda) }{\lambda}  = 0. 
 \end{equation} 
 Typically, for general $H^s$ functions which are not smoother than $H^s$, 
 $\mu \le C \ln \lambda$ with $C$ related to $s$.

 \begin{prop}
 Suppose that $h \in L^2( \RR)$ is not analytic in the sense that it satisfies
 $\eqref{eq511}$. Then  solutions $(u^\lambda, v^\lambda)$ of 
 $\eqref{eq57}$ converge weakly to $(0, h)$.
 \end{prop} 
 
 This means that the weak limits satisfy 
 \begin{equation}
\label{eq512}
\left\{\begin{aligned}
&\partial_t \uu =  0  ,  \qquad \uu_{\vert t = 0} = 0, 
\cr
& \partial_t \uv = 0   ,   \qquad \uv_{\vert t = 0} = h. 
\end{aligned}\right.
\end{equation}
These ``limit'' equations have nothing to see with the original ones \eqref{eq56}, 
implying that \eqref{eq57} are {\sl not} approximations of \eqref{eq56}.

\begin{proof} 
The estimate \eqref{eq510b} and \eqref{eq511} imply that   for all 
$t > 0$, $U^\lambda(t) \to 1$ when 
$\lambda \to \infty$. Since $U^\lambda < 1$ and is increasing, 
this implies that $A^\lambda(t) \to 0 $, uniformly on compacts subsets of $]0, + \infty[$. 
 By \eqref{eq59} 
$$
\widehat u^\lambda(t, \xi) \  \to \  0\qquad  \hat v^\lambda(t, \xi) \to \hat (\xi), 
$$   
uniformly on compacts of $]0, + \infty[\times \RR$. 
 \end{proof}

 \bigbreak
 \begin{remark}
 \textup{The same analysis applies to more general  initial data. 
 One can for instance take  $u(0, \cdot)  = h \ne 0$, with 
 $\Vert h \Vert_{L^2} < 1$ and $v(0, \cdot) = 0$. This only amounts to interchange
 $\cosh$ and $\sinh $ in \eqref{eq59}.}
 \end{remark}


\bibliographystyle{amsalpha}

\begin{thebibliography}{A}

\bibitem{AS} P. D'Ancona, S.Spagnolo, {\it Global solvability for the degenerate
Kirchhoff equation with real analytic data}, Invent. Math., 108 (1992), pp 247-262.  

\bibitem{BGT} S.Baouendi, C.Goulaouic, F.Treves, {\it Uniqueness is certain first order nonlinear complex Cauchy problems}, Comment. Pure Appl. Math., 38 (1985), pp 109-123. 

\bibitem{Be} S.Benzoni-Gavage, {\it Stability of subsonic planar phase boundaries in a 
van der Waals fluid}, Nonlinear Anal., 31(1998), pp 243-263.  



\bibitem{Bo}  J.M.Bony, {\it Calcul symbolique et propagation des singularit\'es pour les 
\'equations aux d\'eriv\'ees partielles non lin\'eaires},  Ann. Sc. Ecole Norm. Sup, 14(1981), 
pp 209-246.

\bibitem{Del} J.M.Delort, {\it F.B.I. Transformation}, Lectures Notes in Math., 1522 (1992), 
Springer Verlag. 

\bibitem{Fan} H.Fan, {\it On a model of the dynamics of liquid/vapor phase transitions}, 
SIAM J.  Appl. Math., 60(2000), pp 1270-1301. 

\bibitem{Fr} H.Freist\"uhler, {\it Some results on the stability of  non-classical   shock waves}, 
J. Partial Diff. Equ., 11(1998), pp 25-38. 

\bibitem{Ha}  J.Hadamard, {\it Lectures on Cauchy's Problem in Linear Partial 
Differential Equations},  Yale Univ. Press, New Haven-London, 1923; {\it Le probl\`eme de Cauchy et les \'equations
aus d\'eriv\'ees partielles hyperboliques}, Herman, Paris, 1932. 

\bibitem{Ha1}  J.Hadamard, {\it  Quelques cas d'impossibilit\'e du probl\`eme de Cauchy}, 
Kazan\v{\i},  In memoriam N.I.Lobatchewsky 2, (1927), pp 163-176. 


\bibitem{Ho}  L.H\"ormander, {\it  The Analysis of Linear Partial Differential Operators}, 
1985, Springer Verlag.
 
 \bibitem{Hou} J.Hounie, J.R.dos Santos Filho {\it  Well-posed Cauchy problems
 for complex nonlinear equations must be semilinear}, Math.Ann., 294 (1992), pp 439-447. 
 
 
 \bibitem{HW} D.Y.Hsieh, X.P.Wang, {\it Phase transitions in van der Waals fluids}, 
 SIAM J. Appl. Math., 57(1997), pp 871-892. 
 
 
 
 \bibitem{IP}  V.Ivrii, V.Petkov,  {\it Necessary conditions for the correctness of the Cauchy problem
 for non-strictly hyperbolic equations},  
 Uspehi Mat. Nauk 29(1974) pp 3-70; also  Russian Math Surveys 29(1974) pp 1-70. .

\bibitem{Jo}  F.John, {\it  Continuous Dependence on Data for Solutions of Partial Differential Equations},  Comm. on Pure and Appl. Math., 13(1960), pp 551-585.

\bibitem{K} G.Kirchhoff, {\it Vorlesungen \"uber Mechanik}, Leipzig, Teubner 1883. 

\bibitem{La}  P.Lax,  {\it Asymptotic solutions of oscillatory initial value problems}, 
Duke Math. J., 24(1957), pp 627-646. 



\bibitem{La1}  P.Lax,  {\it  Non Linear Hyperbolic Equations},   Comm. Pure and Appl. Math., 
6(1953), pp 231-358. 

\bibitem{Li} J.L.Lions, {\it On some questions in boundary value problems of mathematical physics},
in   Contemporary developments in continuum mechanics 
and partial differential equations,  North-Holland, Amsterdam-New York, 1978. 

\bibitem{Me1} G.M\'etivier, {\it  Uniqueness and approximation of solutions if firsto nonlinear 
Equations}, 
Inv.Math., 82(1985), pp 263-282.

\bibitem{Me2} G.M\'etivier, {\it  Counterexamples to Holmgren's uniquess for analytic nonlinear Cauchy problems}, 
Inv.Math., 112(1993), pp 217-222. 

\bibitem{Mi}  S.Mizohata, {\it Some remarks on the Cauchy problem}, 
 J.Math. Kyoto Univ., 1(1961), pp 109-127. 
 
 \bibitem{MN}  C.Morrey,   {\it  On the analyticity of the solutions of analytic non-linear elliptic systems of partial differential equations. II. Analyticity at the boundary}, 
Amer. J. Math. 80(1958),  pp 219-237. 

 \bibitem{Ni1}  T.Nishitani,   {\it   On the Lax-Mizohata theorem in the analytic and Gevrey classes},
 J. Math. Kyoto Univ., 18(1978), pp 509-521. 
 
  \bibitem{Ni2}  T.Nishitani,   {\it A note on the local solvability of the Cauchy problem}, 
  J. Math. Kyoto Univ., 24(1984), pp 281-284. 
  
\bibitem{ST}  M.Sabl\'e-Tougeron, {\it  R\'egularit\'e microlocale pour des probl\`emes aux limites
non lin\'eaires},  Ann. Inst. Fourier,  36(1986), pp 39-82. 
 
 \bibitem{Sj} J.Sj\"ostrand, {\it Singularit\'es analytiques microlocales},  Ast\'erisque, 95(1982), 
 pp 1-166.  
 
 \bibitem{Sl} M.Slemrod, {\it Admissibility criteria for propagating phase boundaries in a 
 van der Waals fluid}, Arch. Rat. Mech. Anal., 81(1983), pp 301-315. 

\bibitem{Wa}  C.Wagschal,  {\it Le probl\`eme de Goursat non lin\'eaire}, 
 J. Math. Pures et Appl., 58(1979), pp309-337.

\bibitem{Wak} Wakabayashi, {\it The Lax-Mizohata theorem for nonlinear Cauchy problems},  Comm. in Partial Diff. Equ., 26 (2001), pp 1367-1384. 

\bibitem{Yag1} K.Yagdjian, {\it A note on Lax-Mizohata theorem for quasilinear equations},
Comm. in Part.Diff.Equ., 23 (1998),  pp 1111-1122. 

\bibitem{Yag} K.Yagdjian, {\it The Lax-Mizohata theorem for nonlinear gauge invariant 
equations}, Nonlinear Anal. T.MA., 49 (2002), pp 159-175. 
\end{thebibliography}

\end{document}